\documentstyle{article}
\title{Hilbert schemes of a surface and Euler characteristics}
\author{Mark Andrea A.  de Cataldo\thanks{
Partially supported by N.S.F. Grant DMS 9701779 and by an
A.M.S. Centennial Fellowship.}}

\date{September 22, 1998}

\newtheorem{tm}{Theorem}[section]
\newtheorem{lm}[tm]{Lemma}

\newtheorem{rmk}[tm]{Remark}

\newtheorem{??}[tm]{Question}

\font\tenmsb=msbm10
\font\sevenmsb=msbm7
\font\fivemsb=msbm5

\newfam\msbfam
\textfont\msbfam=\tenmsb
\scriptfont\msbfam=\sevenmsb
\scriptscriptfont\msbfam=\fivemsb
\def\Bbb#1{{\fam\msbfam #1}}

\font\teneufm=eufm10
\font\seveneufm=eufm7
\font\fiveeufm=eufm5
\newfam\eufmfam
\textfont\eufmfam=\teneufm
\scriptfont\eufmfam=\seveneufm
\scriptscriptfont\eufmfam=\fiveeufm

\newcommand\ci{\cite}

\newcommand\rat{{\Bbb Q}}
\newcommand\comp{{\Bbb C}}

\newcommand\zed{{\Bbb Z}}
\newcommand\nat{{\Bbb N}}

\newcommand\blacksquare{{\hspace*{\fill} $\Box$}}

\newcommand{\SS}{\Sigma}

\begin{document}
\maketitle
\begin{abstract}
{We use basic algebraic topology and Ellingsrud-Stromme results on the Betti numbers of punctual Hilbert schemes of surfaces to compute a generating
function for the  Euler characteristic numbers of the Douady spaces of ``n-points"  associated with a
complex surface. The projective  case was first proved by L. G\"ottsche.}
\end{abstract}

\section{Introduction}
This note is dedicated to the  determination of the Euler characteristic numbers $e(X^{[n]})$ of the Douady spaces of ``n-points" $X^{[n]}$ associated with a
complex surface $X$; see Theorem \ref{tmen}.

For $X$ algebraic (the Douady space is replaced by the Hilbert scheme) see \ci{go}, \ci{go-so} and \ci{ch}.
In the projective case L. G\"ottsche \ci{go} used P. Deligne's solution to the Weil Conjectures 
and  Ellingsrud-Stromme's  determination of the Betti numbers of 
punctual Hilbert schemes \ci{e-s}. Later, G\"ottsche-Soergel \ci{go-so}
proved the quasi-projective case    using
the methods of perverse sheaves.
 J. Cheah
\ci{ch}  obtained similar results using the theory of 
Mixed Hodge structures and \ci{e-s}.

\medskip
The proof presented here, 
valid without the algebraicity restriction on $X$, uses \ci{e-s} 
and elementary algebraic topology.

\bigskip
{\bf Acknowledgments.} I wish to the thank the Max-Planck-Institut
f\"ur Mathematik in Bonn for its warm hospitality and stimulating environment.
I would like to thank Pavlos Tzermias for useful conversations.

\section{Notation and preliminaries}
The topological Euler characteristic of a topological space $Y$, with respect
to singular cohomology with compact supports, is denoted by $e(Y)$. If a statement involves this quantity $e(Y)$, it is understood
that part of the statement consists of the assertion  that $e(Y)$
is defined.
 
 The term {\em complex surface} denotes, here,  a smooth, connected
and complex analytic surface with  separable topology.

Let $X$ be a complex surface and $n\in \nat^*$.
Denote by  $X^{[n]}$ 
the Douady Space \ci{do} of zero-dimensional analytic subspaces of $X$ of length $n$; in the algebraic case this notion is replaced by the one of 
the Hilbert scheme; note that the ``local differential geometries" 
of the Hilbert scheme and of the Douady spaces run parallel so that
one can translate most of the proofs in \ci{fo} in the analytic context.
The spaces $X^{[n]}$ are connected, complex manifolds of dimension $2n$.
The corresponding Barlet space (cf. \ci{ba}) of zero-dimensional cycles in $X$ of total multiplicity $n$
can be identified with the symmetric product
$X^{(n)} =X^n/\SS_n$,
where  $\SS_n$ denotes the symmetric group over $n$ elements and it acts
by permuting the factors of $X^n$.
There is a natural proper holomorphic map $\pi  : X^{[n]} \to X^{(n)}$ which is a resolution of the singularities of $X^{(n)}$ (cf. \ci{ba}); we call this morphism the Douady-Barlet morphism.

There is a natural stratification for this morphism given by the partitions
of $n$. 
Let $P(n)$ denote the set of partitions of the natural number $n$:
$$
P(n): =\left\{ \nu =(\nu_1, \ldots , \nu_k)\, \left| \right. \; k \in \nat^*, \, \nu_j \in \nat^* , \; \nu_1\geq \nu_2 \geq \ldots  \geq \nu_k >0, \, \sum_{j=1}^k \nu_j=n \right\}.
$$
The cardinality of $P(n)$ is denoted by $p(n)$.
Every  partition $\nu=(\nu_1, \ldots, \nu_k)$ can also be represented by
 a sequence
of $n$ nonnegative integers:
$$
\nu=(\nu_1, \ldots , \nu_k) \Longleftrightarrow \left\{ \alpha_1, \ldots, \alpha_n\right\},
$$
where
$\alpha_i$ is the number of times that the integer $i$ appears in the sequence
$(\nu_1, \ldots, \nu_k)$. Clearly $n=\sum_{i=1}^n{ i\,\alpha_i}$.
The number $k$ associated with $\nu$ is called {\em the length of $\nu$} and it is usually denoted by
$\lambda (\nu)$.
The set $P(n)$ is in natural one to one correspondence with a set
of mutually disjoint subsets of $X^{(n)}$ which are locally closed for the Zariski topology on 
analytic spaces: associate with every $\nu \in P(n)$ the space
$$
S^n_{\nu}X=S^n_{(\nu_1,\ldots, \nu_k)}X:=
\left\{\sum_{j=1}^k{\nu_j x_j} \in X^{(n)} \, | \; x_{j_1}\not= x_{j_2} \; if 
\;
j_1\not= j_2 \right\}.
$$
The space $X^{(n)}$ is the disjoint union of  these locally closed subsets: 
\begin{equation}
\label{strs}
{
X^{(n)} =  {\coprod}_{\nu \in P(n)} S^n_{\nu} X.
}
\end{equation}
Define $\widetilde{S^n_{\nu} X}: = \pi^{-1} (S^n_{\nu} X)$.
The space $X^{[n]}$ is the disjoint union of  these locally closed subsets: 
\begin{equation}
\label{strd}
{
X^{[n]} =  {\coprod}_{\nu \in P(n)} \widetilde{S^n_{\nu} X}.
}
\end{equation}

Given $m \in \nat^*$,
the analytic isomorphism class of the fiber over a point $x\in X=
S^m_{(m)} X$ of the Douady-Barlet 
morphism  $\pi: X^{[m]} \to
X^{(m)}$ are independent of the surface $X$ and of the point 
$x\in X$;  since they are isomorphic to the
fiber in question when $X$ is the affine plane ${\Bbb A}^2=\comp^2$
and the point $x$  is the origin, we denote them
by ${\Bbb A}^{2^{[m]}}_{[0]}$.
\begin{tm}
\label{cell} 
(Cf. \ci{e-s})
$e \left( {\Bbb A}^{2^{[m]}}_{[0]} \right)=p(m)$.
\end{tm}
Results of Fogarty's \ci{fo} (see also \ci{go}), adapted in this analytic context,
ensure that the restriction of $\pi$, $\pi:\widetilde{S^n_{\nu} X}
\to S^n_{\nu}X$, is a locally trivial fiber bundle in the analytic
Zariski topology
with fiber $F_{\nu}$ 
\begin{equation}
\label{fib}
F_{\nu}\simeq\prod_{j=1}^k {\Bbb A}^{2^{[\nu_j]}}_{[0]}.
\end{equation}

We shall deal only  with complex analytic spaces $Z$ for which the 
singular cohomology with compact supports
satisifes $\dim_{\rat} H^*_c(W,\rat) < \infty$; in particular the Euler numbers
$e(Z)$ are defined.
Alexander-Spanier cohomology, singular cohomology 
and sheaf cohomology with constant coefficients 
(all with compact supports)
are naturally isomorphic to each other on
 the spaces that we will consider.

Alexander-Spanier's Theory has good properties with respect to closed subsets.
We shall need two of these properties. In what follows we assume 
that $W, Y, Z$ and $Z_i$ are analytic spaces with finite dimensional
rational cohomology with compact supports.
Let $Z = \coprod_{i=1}^n Z_i$ where  $\{Z_i\}$ is a finite collection
of mutually disjoint locally closed subsets of $Z$.
Then
\begin{equation}
\label{a'}
{
e(Z) = \sum_{i=1}^n e\left(Z_i\right).
}
\end{equation}
Let $f:W \to Z$ be  a  continuous map which is locally trivial,
with fiber $Y$,
with respect to the stratification $\Sigma_i$, i.e. $f^{-1} (Z_i)
\simeq \Sigma_i \times Y$ for every $i=1, \ldots , n$. Then
\begin{equation}
\label{b'}
{
e(W) = e(Y) \, e(Z).
}
\end{equation}

It should be now clear what we should be  aiming at:
{\em reduce the  computation of $e \left( X^{[n]} \right) $ to the  computation
of  $e\left( \widetilde{ S^n_{\nu} X } \right) $
via {\rm (\ref{a'})}.} The numbers 
$e\left (\widetilde{S^n_{\nu} X} \right)$
 will be computed using (\ref{b'}) and (\ref{fib}).
The resulting sum will turn out to be a multinomial expansion.

\section{Statement of the result and the numbers  
$e\left( S^{(n)}_{\nu} X \right)$ }
The following is the main result of this note.
\begin{tm}
\label{tmen}
Let $X$ be a smooth complex surface such that
$\dim_{\rat}H^*_c(X, \rat )<\infty$  and with topological Euler
characteristic $e(X)$. Let  $X^{[n]}$ be the associated Douady spaces
of  ``$n$-points."  The following
 is a generating function for $e\left(X^{[n]}\right)$:
\begin{equation}
\label{en1}
\sum_{n=0}^{\infty}{e\left(X^{[n]}\right)\,q^n} = \prod_{k=1}^{\infty} \left( \frac{1}{1-q^k}
\right)^{e(X)}.
\end{equation}
\end{tm}
\begin{rmk}
{\rm The assumptions of the theorem are satisfied for  any {\em algebraic}
surface.}
\end{rmk}
The proof will be preceded by some elementary calculations.

\bigskip
\noindent
{\bf The structure of $S^n_{\nu}X$.}
Let $\nu=(1, \ldots, 1) \in P(n)$. The stratum $S^n_{(1,\ldots, 1)}X$ is the unique
Zariski open subset in $X^{(n)}$  belonging to the stratification
(\ref{strs}). Let
$D^n$ be the big diagonal in the cartesian product $X^n$. The stratum $S^n_{(1, \ldots, 1)}X$  is the quotient
of $X^n \setminus D^n$ under the natural free action of the symmetric group over $n$ elements $\SS_n$. It follows that 
$S^n_{(1, \ldots, 1)}X$ is smooth of complex dimension $2n$.
Each stratum   $S^n_{\nu}X$, $\nu \in P(n)$,  is built from these basic 
strata $S^{m}_{(1, \ldots, 1)}X$ for $m\leq n$. Here is how.
Let
$\nu=(\nu_1, \ldots, \nu_k) \leftrightarrow\{\alpha_1, \ldots, \alpha_n\}$ be
the same partition in the two different pieces of notation.
There is a natural isomorphism:
\begin{equation}
\label{snnx}
S^n_{\nu} X \simeq \left( 
\prod_{i=1, \, \alpha_i \ne 0}^n \left(S^{\alpha_i}_{(1,\ldots, 1)}X \right)\, \right) 
\setminus \Delta, 
\end{equation}
where $\Delta$ is the closed set 
$$
\Delta:= \left\{(C_1, \ldots, C_l) \in  \prod_{i=1, \, \alpha_i\ne 0}^n S^{\alpha_i}_{(1,\ldots, 1)}X \: |
\: Supp (C_h)\cap Supp (C_j) \not= \emptyset  \:for\: some\: h\not= j \right\}.
$$
It follows that each stratum $S^n_{\nu} X$ is smooth,
connected and of  complex dimension $2k=2\lambda (\nu)$.

\bigskip
\noindent
{\bf The numbers $e(S^n_{\nu}X)$.}
Recalling the convention $0!=1$, we have the following elementary 
\begin{lm}
\label{epart}
Let $X$ be a smooth complex analytic surface with
$\dim_{\rat} H^*(X, \rat)  < \infty$, $e=e(X)$, $n\in \nat^*$ and $P(n) \ni \nu=(\nu_1, \ldots, \nu_k)\leftrightarrow
\{\alpha_1, \ldots, \alpha_n\}$. Then
\begin{equation}
\label{esnnx}
e \left( S^n_{\nu} X \right) =   \frac{1}{\alpha_1 ! \ldots \alpha_n  !}  \: e\,(e-1)\, \ldots 
\,(e- (\lambda (\nu ) -1)). 
\end{equation}
\end{lm}
{\em Proof.} 
We first compute $e \left( S^m_{(1,\ldots ,1)} X \right)$ for every $m$. 
Since  $S^m_{(1,\ldots ,1)} X = \left(X^m\setminus D^m \right)/ \SS_m $,
and the action is free, we have that, by the multiplicativity
of the number $e$ under finite covering maps, 
$m! \;e \left( S^m_{1,\ldots ,1} X \right)= e \left(X^m\setminus D^m \right)$.

\smallskip
\noindent
We now prove that  $e \left( X^m\setminus D^m \right)=$
$e(e-1) \ldots (e - (m-1)$ by
induction on $m$. The case $m=1$ is trivial.
Assume that the above formula is true for $m-1$ and let us prove it for $m$.
We have a commutative  diagram 
$$\begin{array}{lcc}
X^m\setminus D^m  & \hookrightarrow & (\, X^{m-1}\setminus D^{m-1}\,) \times X \\
 \\
\quad \quad\downarrow proj_{1,\ldots, m-1}  & \quad \swarrow proj_1 & \,\\
\\
X^{m-1}\setminus D^{m-1} & \,  & \, 
\end{array}
$$
where the horizontal map is the natural open embedding, the vertical map
is the projection to the first $(m-1)$ coordinates and the ``south-west"
map is the projection onto the first factor. The map $proj_1$
admits $(m-1)$ natural sections $(x_1, \ldots, x_{m-1}) \to$
$(x_1, \ldots, x_{m-1}, x_l)$, where $l=1, \ldots , m-1$.
The images of these sections are pairwise disjoint. Their union is a closed subset of $(X^{m-1}\setminus D^{m-1}) \times X$ isomorphic to $(m-1)$
disjoint copies of $X^{m-1}\setminus D^{m-1}$; the complement of this set
is the open subset which we identify with  $X^m\setminus D^m$ via the open immersion in the diagram above.
The additive and multiplicative properties of the number $e$ and the inductive hypothesis
give that 
\begin{eqnarray}
\label{c}
e \left(X^m\setminus D^m \right) & = &
\left[ \, e(e-1) \ldots (e-(m-1-1)) e\, \right] - \left[ \, e(e-1) \ldots (e-(m-1-1)) (m-1) \, \right]  \\
\, & = &
 e(e-1) \ldots (e-(m-1)).
\end{eqnarray} 
This proves that
\begin{equation}
\label{ann}
m! \;e \left( S^m_{1,\ldots ,1} X \right)= e \left(X^m\setminus D^m \right)=
e(e-1) \ldots (e-(m-1)).
\end{equation}
As to the case of $S^n_{\nu} X$ we note that, because of the isomorphism
given in (\ref{snnx}), this space is the quotient
of 
$X^{ \alpha_1 + \ldots +\alpha_n= \lambda (\nu) } \setminus D^{\lambda (\nu)}$ under
the  free action of the group $\SS_{\alpha_1} \times \ldots \times \SS_{\alpha_n}$ induced by the individual actions of the groups
$\SS_{\alpha_i}$ on the factors $X^{\alpha_i}$. The formula follows from the multiplicativity
of Euler characteristics under finite covering maps. \blacksquare

\begin{rmk}
{\rm 
It is an amusing exercise to derive Macdonald formula for the Euler characteristics of symmetric products using Lemma
\ref{epart} 
and the multinomial expansion; i.e. prove that
\begin{equation}
\label{mac}
\sum_{n=0}^{\infty} e\left(X^{(n)}\right) \, q^n =\left( \frac{1}{1-q} \right)^{e(X)}.
\end{equation}
}
 \end{rmk}

\bigskip
\noindent
{\bf The numbers $e(\widetilde{S^n_{\nu}X})$.}
By virtue of Lemma \ref{epart}, Theorem  \ref{cell} and (\ref{fib}), 
the basic additive and multiplicative properties 
of cohomology with compact support
give the following
\begin{lm}
\label{etilde} Let $X$ be a complex surface
with finite Betti numbers. Then
\begin{eqnarray}
e \left(
\widetilde{S^n_{\nu} X} \right) & = & e\left(\prod_{j=1}^{\lambda (\nu) } {\Bbb A}^{2^{[\nu_j]}}_{[0]}\right)\:
e(S^n_{\nu} X)=  p(\nu_1) \ldots  p(\nu_k) \;\; 
\frac{1}{\alpha_1 ! \ldots \alpha_n  !}  \: e\,(e-1)\, \ldots 
\,(e- (\lambda (\nu ) -1))  \nonumber\\
& = &
p(1)^{\alpha_1} \ldots  p(n)^{\alpha_n} \;\; 
\frac{1}{\alpha_1 ! \ldots \alpha_n  !}  \: e\,(e-1)\, \ldots 
\,(e- (\lambda (\nu ) -1)) \, .
\end{eqnarray}
\end{lm}

\section{Proof of Theorem \ref{tmen}}
{\em Proof of Theorem \ref{tmen}.}
By virtue of  (\ref{strd}), (\ref{a'}) and Lemma
\ref{etilde} we have
\begin{eqnarray}
e\left( X^{[n]}\right) & = & \sum_{\nu \in P(n)}  e\left( \widetilde{ S^n_{\nu} X } \right)  \\
& = &
\sum_{k=1}^n{
\sum_{ \begin{array}{c} 
\nu \in P(n) \\
\lambda (\nu )=k
\end{array}
} 
\frac{p(1)^{\alpha_1} \ldots \,  p(n)^{\alpha_n}}{\alpha_1  ! \ldots 
\alpha_n !} \: e\,(e-1)\, \ldots 
\,(e- (k -1))
}.
\end{eqnarray}
The last summand is the expansion of
the r.h.s.  of equation (\ref{en1}). Let us check this for $e\geq 0$
from a combinatorial point of view. Set $p(0):=1$.
For $e\geq 0$ the product $\prod_{k=1}^{\infty} \left(\frac{1}{1-q^k}\right)^{e(X)}$
can be  defined as the limit, in $\rat [[q]]$, of the polynomials
$P_n(q):=(\sum_{s=0}^n{ p(s)\,q^s})^e$.
The coefficients $c_{m,s}$ of $q^s$ in $P_{m}(q)$ are the same as the one of the infinite product for every $s\leq n$ and every $m\geq n$.
For $n\leq e$ their form is given by the multinomial coefficient expansion
in the following way:
\begin{eqnarray}
c_{m,s} & = &  \sum_{
\begin{array}{c}
\sum_{i=0}^n a_i= e \\
 \sum_{i=0}^n i\,a_i= n
\end{array}  
}{e\choose {a_0, \ldots , a_n}} p(0)^{a_0}
\ldots p(n)^{a_n}   \\
 & = &  \sum_{k=1}^n{
\sum_{
\begin{array}{c}
\sum_{j=1}^n a_j= e-k \\
 \sum_{j=1}^n j\,a_j= n 
\end{array}
}
}
{e\choose {k, a_1, \ldots , a_n}}  p(1)^{a_1}
\ldots p(n)^{a_n}  \\
& = &
\sum_{k=1}^n{
\sum_{
\begin{array}{c}
\sum_{j=1}^n a_j= e-k \\
 \sum_{j=1}^n j\,a_j= n 
\end{array}
}
}
\frac{p(1)^{\alpha_1} \ldots p(n)^{\alpha_n}}{\alpha_1  ! \ldots 
\alpha_n !} \: e\,(e-1)\, \ldots 
\,(e- (k -1)), 
\end{eqnarray}
which coincides with $e\left(X^{[n]}\right)$.
The formula is still correct for $n>e$ (many summands are zero). 
The case $e=0$ is trivial. We leave the case $e<0$ to the reader.
\blacksquare

\begin{rmk}{\rm 
The  generating function
of Theorem \ref{tmen} exhibits modular behavior (cf. \ci{go}).
}
\end{rmk}

\begin{rmk}
{\rm 
If $X$ is not algebraic, then the approach in \ci{ch} cannot be used to prove Theorem \ref{tmen}.  The paper \ci{go-so} 
restricts itself to the algebraic context. However, 
that restriction is unnecessary.
This is remarked in \ci{de-mi}, whose main purpose is to 
prove the necessary decomposition theorem in the form of an explicit quasi-isomorphism
of complexes.}
\end{rmk}

\begin{rmk}
{\rm 
For $X$ algebraic the numbers $e\left(X^{[n]}\right)$ and  the orbifold Euler characteristics
$e(X^n, \SS_n)$ coincide. Similarly, we see that because
of Theorem \ref{tmen} this ``coincidence" occurs for any complex surface, i.e. not necessarily
algebraic. This  is explained by the fact that
the groups $K ( X^{[n]} ) \otimes_{\zed} \rat$ and the $\SS_n$-Equivariant K-Theory
of $X^n$, $K_{\SS_n}(X^n) \otimes_{\zed} \rat$, are {\em naturally} isomorphic.
This fact is due to   B. Bezrukavnikov and  V. Ginzburg \ci{ba-gi} and, independently, 
to L. Migliorini and myself \ci{de-mi}.
}
\end{rmk}

\begin{??} {\rm 
The numbers $e\left(X^{[n]}\right)$ are  invariant in the class of
 complex analytic smooth surfaces homotopically equivalent to $X$.
Is $X^{[n]}$  a topological invariant
of $X$? To be precise,  let  $X$ and $Y$ be  complex surfaces homeomorphic
to each other;  is $X^{[n]}$ homeomorphic to $Y^{[n]}$? Is $X^{[n]}$ an invariant for
the differentiable structure of $X$?
}
\end{??}



\bigskip
\noindent
Author's address:
Department of Mathematics, Harvard University, Science Center, One Oxford Street,
Cambridge, MA 02138, U.S.A. $\quad$
e-mail: {\em mde@math.harvard.edu}

\end{document}